\documentclass[a4paper, amsmath, amssymb, 11pt]{article}

\usepackage[top=4cm, bottom=3.5cm, left=3cm, right=3cm]{geometry}
\usepackage[fleqn]{amsmath}
\usepackage{amssymb}
\usepackage{amsthm}
\usepackage{mathtools}
\begin{document}

\begin{center}
\large {\bf A Generalized Function defined by the Euler first kind integral
and its connection with the Dirac delta function}

\vspace{0.4cm}

\large Vagner Jikia$^{1}$, Ilia Lomidze$^{2}$ \\
{\small $^{1}$ Javakhishvili State University, Tbilisi, Georgia, \\
$^{2}$ Georgian Teknical University, Engineering Physics Department, Tbilisi, Georgia. \\ v\_jikia@yahoo.com, lomiltsu@gmail.com }
	
\end{center}
%
%\begin{abstract}
%
%{\small As it is known, quantum mechanical formalism of two charge particles problem of %continuous spectra is incomplete and inconsistent up to now. In present article we attempt to %remove partially this distortion in Schrodinger formalism frame. There are proposed the main %results of our research connected to the Coulomb quantum mechanical problem.}
%\end{abstract}
%\vspace{0.2cm}

{\bf Abstract.} We have shown that in some region where the Euler integral of the first kind diverges, the Euler formula defines a generalized function. The connected of this generalized function with the Dirac delta function is found. \\



It is known (see, e.g. [1]) that the Euler beta function (the Euler integral of the first kind)
\begin{eqnarray}  
&B(\alpha ,\,\beta ) = \int_0^1 {{t^{\alpha  - 1}}{{(1 - t)}^{\beta  - 1}}dt} , \\
&\nonumber{\mathop{\rm Re}\nolimits} \,(\alpha ) > 0,\,\,{\mathop{\rm Re}\nolimits} \,(\beta ) > 0,
\end{eqnarray}
satisfies the Euler formula:
\begin{eqnarray}  
\,B(\alpha ,\,\beta ) = \frac{{\Gamma (\alpha )\Gamma (\beta )}}{{\Gamma (\alpha  + \beta )}},
\end{eqnarray}
where $ \Gamma (z) $  is the Euler gamma function (the Euler integral of the second kind):
\begin{eqnarray}  
\nonumber\Gamma (z) = \int_0^1 {{t^{z - 1}}\exp ( - t)dt} , ~~ {\mathop{\rm Re}\nolimits} z > 0.
\end{eqnarray}

Based on the known property of the gamma function $ \Gamma (z + 1) = z\Gamma (z), $
%\begin{eqnarray}  
%\Gamma (z + 1) = z\Gamma (z),
%\end{eqnarray}
the formula (2) can be used to derive a recurrent relations $ \left( {n,\,k = 1,2,...} \right) $
\begin{eqnarray}  
&\nonumber B(\alpha  + n,\beta  + k) = \dfrac{{\Gamma (\alpha  + n)\Gamma (\beta  + k)}}{{\Gamma (\alpha )\Gamma (\beta )}}\dfrac{{\Gamma (\alpha  + \beta )}}{{\Gamma (\alpha  + \beta  + n + k)}}B(\alpha ,\beta ) \\
&\nonumber {\mathop{\rm Re}\nolimits} \,(\alpha ) > 0,\,\,{\mathop{\rm Re}\nolimits} \,(\beta ) > 0,
\end{eqnarray}
and
\begin{eqnarray}  
&\nonumber B(\alpha  - n,\beta  - k) = \dfrac{{\Gamma (\alpha  - n)\Gamma (\beta  - k)}}{{\Gamma (\alpha )\Gamma (\beta )}}\dfrac{{\Gamma (\alpha  + \beta )}}{{\Gamma (\alpha  + \beta  - n - k)}}B(\alpha ,\beta ), \\
&\nonumber {\mathop{\rm Re}\nolimits} \,(\alpha ) > n,\,\,{\mathop{\rm Re}\nolimits} \,(\beta ) > k.
\end{eqnarray}
These formulae, as it is well-known, can be used to define the beta function on the whole complex plane of its arguments, starting from the domain $ 0 < {\mathop{\rm Re}\nolimits} \,(\alpha ) \le 1,\;\;0 < {\mathop{\rm Re}\nolimits} \,(\beta ) \le 1, $
except the lines $ {\mathop{\rm Re}\nolimits} \alpha  = 0, - 1, - 2,...,\,\quad {\mathop{\rm Re}\nolimits} \beta  = 0, - 1, - 2,...\, $ 

An analytical continuation of $ B(\alpha ,\,\beta ) $  on the whole complex plane of its arguments which is made by Pochhammer using a double loop integral (see, e.g., [1], \textsection 12.43): 
\begin{eqnarray}  
\exp [ - i\pi (\alpha  + \beta )]\int\limits_P^{(1 + ,0 + ,1 - ,0 - )} {\,{t^{\alpha  - 1}}{{(1 - t)}^{\beta  - 1}}dt}  =  - 4\sin (\pi \alpha )\sin (\pi \beta )B(\alpha ,\,\beta )
\end{eqnarray}
gives us an opportunity to find the divergence of $ B(\alpha ,\,\beta ) $  at the lines ${\mathop{\rm Re}\nolimits} \alpha ,\beta  = 0, - 1, - 2,...\,\,. $

According to definition of the Dirac delta function (see e.g. [2], \textsection 5.7) we assume that it is defined as weak limit of a sequence of some approximate functions $ {\omega _\varepsilon }(x): $
\begin{eqnarray}  
\delta (x)\mathop  = \limits^{weak} \mathop {\lim }\limits_{\varepsilon  \to 0 + } {\omega _\varepsilon }(x),
\end{eqnarray}
where approximate function $ {\omega _\varepsilon }(x) $  can be constructed as (see [2], \textsection 5.2)
\begin{eqnarray}  
{\omega _\varepsilon }(x) = {\varepsilon ^{ - 1}}\eta ({\varepsilon ^{ - 1}}x),
\end{eqnarray}
for any bounded function $ \eta (x), $  such that
\begin{eqnarray}  
\int_{ - \infty }^\infty  {\eta (x)dx}  = 1.
\end{eqnarray} 
In the next consideration we use functions $ \eta (x) = {(\pi x)^{ - 1}}\sin x $  and  $ \eta (x) = {[\pi (1 + {x^2})]^{ - 1}}. $

Another mostly symbolic but rather common definition (see e.g. [2], \textsection 5.1; [3], \textsection 2.6) considers the Dirac delta as generalized function which satisfies the (weak) equality 
\begin{eqnarray}
\int_a^b {\varphi (x)\delta (x)dx}  = \left\{ {\begin{array}{*{20}{c}}
{0, ~~~~~~ 0 \notin [a,b]} \\[0.4em]
{\varphi (0), ~~ 0 \in (a,b)}
\end{array}} \right.\
\end{eqnarray} 
for any continuous function $ \varphi (x) $ defined on the interval $ [a,b]. $ 

There was shown (see [4] and [5]) that the next theorem is fulfilled:

{\bf Theorem.} {\it For $ {\mathop{\rm Re}\nolimits} \alpha  = {\mathop{\rm Re}\nolimits} \beta  = 0 $  and $ {\mathop{\rm Im}\nolimits} \alpha  =  - {\mathop{\rm Im}\nolimits} \beta  = x \in \mathbb{R}, $  the Euler beta function   turns into generalized function and can be expressed through the Dirac delta function:}
\begin{eqnarray}  
B(ix, - ix) = \mathop {\lim }\limits_{\varepsilon  \to 0 + } B(\varepsilon  + ix,\varepsilon  - ix) = 2\pi \delta (x).
\end{eqnarray} 

Let us prove the next lemma before proving the Theorem:

{\bf Lemma.} {\it Let $ f(\varepsilon ,x) $  be a function of two real variables,  
$ \varepsilon  \in (a,b),\ a < 0 < b,\;\,x \in \mathbb{R}, $   which is continuous on the area of definition and has continuous first derivative $ \frac{\partial }{{\partial \varepsilon }}f(\varepsilon ,\;x) = {f'_\varepsilon }(\varepsilon ,\;x) $ and continuous and limited second derivative 
\begin{eqnarray}  
\frac{{{\partial ^2}}}{{\partial {\varepsilon ^2}}}f(\varepsilon ,\;x) = {f''_{\varepsilon \varepsilon }}(\varepsilon ,\;x) \le M < \infty 
\end{eqnarray}
everywhere on this area. Then}
\begin{eqnarray}  
\mathop {\lim }\limits_{\varepsilon  \to 0 + } \left( {f(\varepsilon ,x)\frac{\varepsilon }{{{\varepsilon ^2} + {x^2}}}} \right) = \pi f(0,0)\delta (x).
\end{eqnarray}

{\bf Proof of the lemma.} Expanding the function $ f(\varepsilon ,x) $ in the neighborhood of the point $ (0,x) $ where $ x \in \mathbb{R}, $ according to the Taylor formula with remainder in the Lagrange form
\begin{eqnarray}  
\nonumber f(\varepsilon ,\;x) = f(0,\;x) + \frac{1}{{1!}}{f'_\varepsilon }(0,\;x)\varepsilon  + \frac{1}{{2!}}{f''_{\varepsilon \varepsilon }}(\xi ,\;x){\varepsilon ^2},\quad 0 < \xi  < \varepsilon ,\quad x \in \mathbb{R},
\end{eqnarray}
one obtains
\begin{eqnarray}  
&\nonumber f(\varepsilon ,x)\dfrac{\varepsilon }{{{\varepsilon ^2} + {x^2}}} = \left( {f(0,x) + \dfrac{1}{{1!}}{{f'}_\varepsilon }(0,x)\varepsilon  + \dfrac{1}{{2!}}{{f''}_{\varepsilon \varepsilon }}(\xi ,\;x){\varepsilon ^2}} \right)\pi {\omega _\varepsilon }(x), \\
&\nonumber 0 < \xi  < \varepsilon , \ x \in \mathbb{R}, 
\end{eqnarray}
where the approximate function $ {\omega _\varepsilon }(x) = {(\pi \varepsilon )^{ - 1}}{[1 + {(x/\varepsilon )^2}]^{ - 1}} $ obviously satisfies (5) and (6). Calculating the weak limit $ \varepsilon  \to 0 + $ of the both sides of the formula obtained and using the conditions of the Lemma and the equality (7) one gets
\begin{flalign}
\nonumber &~~\mathop {\lim }\limits_{\varepsilon  \to 0 + } \int_{ - \infty }^\infty  {\varphi (x)f(\varepsilon ,x)\frac{\varepsilon }{{{\varepsilon ^2} + {x^2}}}dx}
\\
\nonumber & ~~~~~~~~~~~~~~~~~~~~~~~~~~~~~~ = \mathop {\lim }\limits_{\varepsilon  \to 0 + } \int_{ - \infty }^\infty  {\varphi (x)\left( {f(0,x)\pi {\omega _\varepsilon }(x) + o(\varepsilon )} \right)dx}  = \pi \varphi (0)f(0,0),
\end{flalign} 
for any continuous function $ \varphi (x). $ \\
{\bf Proof of the Theorem.} Using the formula (2) and the property of the gamma function $ z\Gamma (z) = \Gamma (1 + z) $
%\begin{eqnarray}
%\nonumber z\Gamma (z) = \Gamma (1 + z)
%\end{eqnarray}
one can write down
\begin{eqnarray}
\nonumber B(\varepsilon  + ix,\varepsilon  - ix) = \frac{{\Gamma (\varepsilon  + ix)\Gamma (\varepsilon  - ix)}}{{\Gamma (2\varepsilon )}} = \frac{{\Gamma (\varepsilon  + ix + 1)\Gamma (\varepsilon  - ix + 1)}}{{\Gamma (2\varepsilon  + 1)}}\frac{{2\varepsilon }}{{{\varepsilon ^2} + {x^2}}}.
\end{eqnarray}
Now to calculate the limit in the formula (8) it is sufficient to note that the factor
\begin{eqnarray}
\nonumber \frac{{2\Gamma (\varepsilon  + ix + 1)\Gamma (\varepsilon  - ix + 1)}}{{\Gamma (2\varepsilon  + 1)}}
\end{eqnarray}
satisfies the conditions of the Lemma: the function
\begin{eqnarray}
f(\varepsilon ,\;x) = \frac{{2\Gamma (\varepsilon  + ix + 1)\Gamma (\varepsilon  - ix + 1)}}{{\Gamma (2\varepsilon  + 1)}}
\end{eqnarray}
is analytic in the point  $ \varepsilon  = 0,~x = 0 $ of the complex plane, $ f(\varepsilon ,x) \in {C^\infty } $  and
\begin{eqnarray}
\nonumber& f(0,\;x) = 2\Gamma (ix + 1)\Gamma (1 - ix), & \,\,\,\,\,\,\,\,\,\,\,\,f(0,\;0) = 2,\\
\nonumber& {{f'}_\varepsilon }(0,\;x) = 2f(0,\;x)\left[ {{\mathop{\rm Re}\nolimits} \psi (1 + ix) - \psi (1)} \right], & \,\,\,\,\,\,\,\,\,{{f'}_\varepsilon }(0,\;0) = 0,\\
\nonumber& {{f''}_{\varepsilon \varepsilon }}(\xi ,\;x) = 2f(\xi ,\;x)\left\{ {2{{\left[ {{\mathop{\rm Re}\nolimits} \psi (\xi  + ix + 1) - \psi (2\xi  + 1)} \right]}^2}} \right.\\
\nonumber& ~~ + \left. {\left[ {{\mathop{\rm Re}\nolimits} \psi '(\xi  + ix + 1) - 2\psi '(2\xi  + 1)} \right]} \right\}, & \,\,\,\,\,{{f''}_{\varepsilon \varepsilon }}(\xi ,\;x) < M,
\end{eqnarray}
where $ \psi (z) $ is the digamma function which is determined as [5]
\begin{eqnarray}
\nonumber \psi (z) = \frac{d}{{dz}}\ln \Gamma (z) = \frac{{\Gamma '(z)}}{{\Gamma (z)}},\quad \psi (1) =  - \gamma ,
\end{eqnarray}
and  $ \gamma $ is the Euler-Mascheroni constant.

Therefore, one obtains the statement to be proven (8).

{\bf Corollary 1.} Using the explicit representation (1)
\begin{eqnarray}
\nonumber B(\varepsilon  + ix,\varepsilon  - ix) = \int_0^1 {dt{t^{\varepsilon  + ix - 1}}{{(1 - t)}^{\varepsilon  - ix - 1}}} ,
\end{eqnarray}
and applying the formula (8) one obtains
\begin{eqnarray}
\mathop {\lim }\limits_{\varepsilon  \to 0 + } \int_0^1 {dt{t^{\varepsilon  + ix - 1}}{{(1 - t)}^{\varepsilon  - ix - 1}}}  = 2\pi \delta (x).
\end{eqnarray}

On the other hand, using an insertion
\begin{equation}
~~~~~~~~~~~~~~~~~ t = \frac{{1 - {\mathop{\rm th}\nolimits} (\xi /2)}}{2},\quad   \xi  = \ln \dfrac{{1-t}}{t},\quad  d\xi  =  - \frac{{dt}}{{t(1 - t)}}\tag{*},
\end{equation}
it is easy to check that (compare with Note 1)
\begin{eqnarray}
\nonumber \mathop {\lim }\limits_{\lambda ,\mu  \to 0 + } \int_\mu ^{1 - \lambda } {dt{t^{ix - 1}}{{(1 - t)}^{ - ix - 1}}}  = \mathop {\lim }\limits_{\alpha ,\beta  \to  + \infty } \,\,\int_{ - \alpha }^\beta  {\exp (i\xi x)d\xi }  = 2\pi \delta (x),
\end{eqnarray}
and therefore the expression $ \int_0^1{dt{t^{ix - 1}}{{(1 - t)}^{ - ix - 1}}} \, $ determines a
generalized function ([2], \textsection 9.2)
\begin{eqnarray}
\int_0^1 {dt{t^{ix - 1}}{{(1 - t)}^{ - ix - 1}}}  = 2\pi \delta (x)\,.
\end{eqnarray}
Comparing this result with (12) one obtains next

{\bf Corollary 2.} Procedures of the limit and integration in the formula (12) can be permuted:
\begin{eqnarray}
\mathop {\lim }\limits_{\varepsilon  \to 0 + } \,\int_0^1 {{t^{\varepsilon  + ix - 1}}{{(1 - t)}^{\varepsilon  - ix - 1}}dt}  = \int_0^1 {\mathop {\lim }\limits_{\varepsilon  \to 0 + } {t^{\varepsilon  + ix - 1}}{{(1 - t)}^{\varepsilon  - ix - 1}}dt}  = 2\pi \delta (x).
\end{eqnarray}
More correctly, according to ([2], \textsection 5.5), the formula should be
\begin{flalign}
\nonumber & \mathop {\lim }\limits_{\varepsilon  \to 0 + } \,\int_{ - A}^B {\varphi (x)dx} \int_0^1 {{t^{\varepsilon  + ix - 1}}{{(1 - t)}^{\varepsilon  - ix - 1}}dt}
\\
\nonumber & ~~~~~~~~~~~~~~~~~~~~ = \mathop {\lim }\limits_{\lambda ,\mu  \to 0 + } \int_{ - A}^B {\varphi (x)dx} \int_\mu ^{1 - \lambda } {\mathop {\lim }\limits_{\varepsilon  \to 0 + } \,{t^{\varepsilon  + ix - 1}}{{(1 - t)}^{\varepsilon  - ix - 1}}dt}  = 2\pi \varphi (0).
\end{flalign}

{\bf Corollary 3.} According to (8) and taking into account that the Dirac delta is even function one obtains the symmetry property:
\begin{eqnarray}
\nonumber B(ix, - ix) = B( - ix,ix).
\end{eqnarray}
Besides, one can generalize formula (8):
\begin{eqnarray}
B(ix - n, - ix - k) = B( - ix - k,ix - n) = \frac{{(n + k)!}}{{n!k!}}2\pi \delta (x),
\end{eqnarray}
where $ k $ and $ n $ are any nonnegative integers. 

{\bf Proof of the Corollary 3.} According to the formulae (2) and (11) one has
\begin{eqnarray}
\begin{array}{l}
\nonumber B(\alpha  - n,\beta  - k) = \dfrac{{\Gamma (\alpha  - n)\Gamma (\beta  - k)}}{{\Gamma (\alpha  + \beta  - n - k)}} = \\ [0.8em]
 = \dfrac{{\Gamma (1 + n + k - \alpha  - \beta )}}{{\Gamma (1 + n - \alpha )\Gamma (1 + k - \beta )}}\dfrac{\pi }{{\sin \pi \alpha }}\dfrac{\pi }{{\sin \pi \beta }}\dfrac{{\sin \pi (\alpha  + \beta )}}{\pi } = \\ [0.8em]
 = \dfrac{{\Gamma (1 + n + k - \alpha  - \beta )}}{{\Gamma (1 + n - \alpha )\Gamma (1 + k - \beta )}}\dfrac{{\Gamma (\alpha )\Gamma (1 - \alpha )\Gamma (\beta )\Gamma (1 - \beta )}}{{\Gamma (\alpha  + \beta )\Gamma (1 - \alpha  - \beta )}} 
\\ [0.8em]
 = B(\alpha ,\beta )\dfrac{{\Gamma (1 + n + k - \alpha  - \beta )}}{{\Gamma (1 + n - \alpha )\Gamma (1 + k - \beta )}}\dfrac{{\Gamma (1 - \alpha )\Gamma (1 - \beta )}}{{\Gamma (1 - \alpha  - \beta )}}, 
\end{array}
\end{eqnarray}
for any complex $ \alpha ,\beta $ such that $ {\mathop{\rm Re}\nolimits} \alpha ,\beta  \in (0,1), $ and for any integers $ n, k $ (see the Pochhammer result (3)). Assuming here
\begin{eqnarray}
\alpha  = \varepsilon  + ix, & \beta  = \varepsilon  - ix, & \varepsilon  > 0,
\end{eqnarray}
we get 
\begin{flalign}
\nonumber & B(\varepsilon  + ix - n,\varepsilon  - ix - k) =
\\
\nonumber & ~~~~~~~~~~~~~~~~~~~~~~~~ = \frac{{\Gamma (1 + n + k + 2\varepsilon )\Gamma (1 - \varepsilon  - ix)\Gamma (1 - \varepsilon  + ix)}}{{\Gamma (1 + n - \varepsilon  - ix)\Gamma (1 + k - \varepsilon  + ix)\Gamma (1 - 2\varepsilon )}}B(\varepsilon  + ix,\varepsilon  - ix),
\end{flalign}
and then, using the techniques similar to the one above to calculate the limit $ \varepsilon  \to 0 + , $ we obtain the equality for generalized functions (see, e.g. [2], \textsection 5.5)
\begin{flalign}
\nonumber & B(ix - n, - ix - k) = \mathop {\lim }\limits_{\varepsilon  \to 0 + } \,B(\varepsilon  + ix - n,\varepsilon  - ix - k)=
\\
\nonumber & ~~~~~~~~~~~~~~~~~~~~~~~~ = \frac{{\Gamma (1 + n + k)}}{{\Gamma (1 + n - ix)\Gamma (1 + k + ix)}}\frac{{\Gamma (1 - ix)\Gamma (1 + ix)}}{{\Gamma (1)}}B(ix, - ix).\,
\end{flalign} 
Taking here into account the formula (8) we finally obtain (15).

	Note, that the right hand side of the formula (15) is not analytical function, therefore it can not be generalized for nonintegers $ n, k. $
	
	Let $ \alpha ,\beta ,x,y \in \mathbb{R}, $ and let consider 
\begin{eqnarray}
\mathop {\lim }\limits_{\alpha ,\beta  \to 0 + } B(\alpha  + ix,\beta  + iy) = \mathop {\lim }\limits_{\alpha ,\beta  \to 0 + } \int_0^1 {{t^{\alpha  + ix - 1}}{{(1 - t)}^{\beta  + iy - 1}}dt} .
\end{eqnarray}
Using here the Euler formula (2) one obtains
\begin{flalign}
\nonumber & B(\alpha  + ix,\beta  + iy) = \frac{{\Gamma (\alpha  + ix)\Gamma (\beta  + iy)}}{{\Gamma [\alpha  + \beta  + i(x + y)]}} =
\\
& ~~~~~~~~~~~~~~~~~~~~~~~~ = \frac{{\Gamma (\alpha  + ix + 1)\Gamma (\beta  + iy + 1)}}{{\Gamma [\alpha  + \beta  + i(x + y) + 1]}} \cdot \frac{{\alpha  + \beta  + i(x + y)}}{{(\alpha  + ix)(\beta  + iy)}}.
\end{flalign}
The first factor in the right hand side of the formula (17) is regular and equal to 1, if $ \alpha  + ix = 0, \beta  + iy = 0; $  it is easy to check that it satisfies the conditions similar to the conditions of the Lemma respect to all (real) variables $ \alpha ,x,\beta ,y. $ 

And what about the second factor,
\begin{eqnarray}
\nonumber \frac{{\alpha  + \beta  + i(x + y)}}{{(\alpha  + ix)(\beta  + iy)}} = \frac{1}{{\alpha  + ix}} + \frac{1}{{\beta  + iy}},
\end{eqnarray}
according to the Sokhotski formula (see, e.g. [2], \textsection 5.8), one has
\begin{flalign}
\nonumber & \mathop {\lim }\limits_{\alpha ,\beta  \to 0 + } \frac{{\alpha  + \beta  + i(x + y)}}{{\alpha \beta  - xy + i(\beta x + \alpha y)}} =  - i\mathop {\lim }\limits_{\alpha ,\beta  \to 0 + } \left[ {\frac{1}{{x - i\alpha }} + \frac{1}{{y - i\beta }}} \right]
\\
& ~~~~~~~~~~~~~~~~~~~~~~~~~~~~~~~~~~~~~~~ = \pi \left[ {\delta (x) + \delta (y)} \right] - i\left[ {P\left( {\frac{1}{x}} \right) + P\left( {\frac{1}{y}} \right)} \right].
\end{flalign}
Hence, calculating the limit in the expression (16), finally we obtain
\begin{flalign}
\nonumber & \mathop {\lim }\limits_{\alpha ,\beta  \to 0 + } B(\alpha  + ix,\beta  + iy) = \pi \left[ {\delta (x) + \delta (y)} \right]
\\
& ~~~~~~~~~~~~~~~~~~~~~~~~~~~~~~~~~ - i\frac{{\Gamma (1 + ix)\Gamma (1 + iy)}}{{\Gamma [1 + i(x + y)]}}\left[ {P\left( {\frac{1}{x}} \right) + P\left( {\frac{1}{y}} \right)} \right], \\
\nonumber & ~~~~~~~~~~~~~~~~~~~~~~~~~~~~~~~~~~ \left( {x + y \ne 0} \right).
\end{flalign}
In the case $ x + y = 0 $ the variables $ x $ and $ y $ are not independent. Therefore this case has to be considered separately. One gets
\begin{flalign}
\nonumber & \mathop {\lim }\limits_{\alpha ,\beta  \to 0 + } {\frac{{\alpha  + \beta  + i(x + y)}}{{\alpha \beta  - xy + i(\beta x + \alpha y)}}_{|x = - y}} =
\\
\nonumber & ~~~~~~~~~~~~~~~~~~~ = \mathop {\lim }\limits_{\alpha ,\beta  \to 0 + } \frac{{\alpha  + \beta }}{{\alpha \beta  + {x^2} - i(\alpha  - \beta )x}} =  - i\mathop {\lim }\limits_{\alpha ,\beta  \to 0 + } \left[ {\frac{1}{{x - i\alpha }} + \frac{1}{{ - x - i\beta }}} \right],
\end{flalign}
and, using the Sokhotski formula again, one obtains:
\begin{eqnarray}
\mathop {\lim }\limits_{\alpha ,\beta  \to 0 + } B(\alpha  + ix,\beta  - ix) = 2\pi \delta (x),
\end{eqnarray}
regardless of $ \alpha  = \beta  = \varepsilon  \to 0, $ or the parameters $ \alpha $ and $ \beta $ tend to zero independently.

On the other hand, using the insertion $ \left( * \right) $, as like as above one obtains (see Note 2)
\begin{eqnarray}
\nonumber & B(ix,iy) \equiv \int_0^1 {{t^{ix - 1}}{{(1 - t)}^{iy - 1}}dt}  = \int_0^1 {{{\left( {\frac{{1 - t}}{t}} \right)}^{i(y - x)/2}}{{[t(1 - t)]}^{ - 1 + i(y + x)/2}}dt} = \\
\nonumber & = \int\limits_{ - \infty }^\infty  {\exp [i\xi (y - x)/2]{{[2\cosh (\xi /2)]}^{i(y + x)}}d\xi ,}
\end{eqnarray}
and (13) follows from this again. Thus, the generalization of the permutation property (14) is valid:
\begin{eqnarray}
\mathop {\lim }\limits_{\alpha ,\beta  \to 0 + } \,\,\int_0^1 {{t^{\alpha  + ix - 1}}{{(1 - t)}^{\beta  - ix - 1}}dt}  = \int_0^1 {\mathop {\lim }\limits_{\alpha ,\beta  \to 0 + } \,\,{t^{\alpha  + ix - 1}}{{(1 - t)}^{\beta  - ix - 1}}dt\,}  = 2\pi \delta (x).
\end{eqnarray}

A number of rather interesting relations which follow to the results obtained above will be considered in detail in the forthcoming publications. \\

{\bf Acknowledgement} This work was supported in part by Georgian Shota Rustaveli National Science Foundation (grant FR/417/6-100/14)  
\\

{\bf Note 1.} One has:
\begin{flalign}
\nonumber & \int_0^1 {dt{t^{ix - 1}}{{(1 - t)}^{ - ix - 1}} = } \mathop {\lim }\limits_{\scriptstyle\lambda  \to 0 + \hfill\atop
\scriptstyle\mu  \to 0 + \hfill} \int_\mu ^{1 - \lambda } {dt{t^{ix - 1}}{{(1 - t)}^{ - ix - 1}}}  = \mathop {\lim }\limits_{\scriptstyle\alpha  \to  + \infty \hfill\atop
\scriptstyle\beta  \to  + \infty \hfill} \int_{ - \alpha }^\beta  {\exp (i\xi x)d\xi } .
\\
\nonumber & ~~~~ \int_{ - \alpha }^\beta  {\exp (i\xi x)d\xi }  = \frac{1}{{ix}}\left( {\exp (i\beta x) - \exp ( - i\alpha x)} \right) = \\
\nonumber & ~~~~~~~~~~~~~~~~~~~~~~~~~~~~~~~~~~  = \frac{1}{{ix}}\left( {\cos (\beta x) - \cos (\alpha x)} \right) + \frac{1}{x}\sin (\beta x) + \frac{1}{x}\sin (\alpha x).
\end{flalign}
{\bf Statement.} { \it For any continuous function $ \varphi (x) $ one has}
\begin{eqnarray}
\nonumber \int_A^B {\varphi (x)\int_0^1 {dt{t^{ix - 1}}{{(1 - t)}^{ - ix - 1}}} dx} = \left\{ {\begin{array}{*{20}{c}}
{2\pi \varphi (0), \ x = 0 \in (A,B)} \\ [0.4em]
{0, \ x = 0 \notin (A,B)}
\end{array}} \right.
\end{eqnarray}
{\bf Proof.} \\
For any continuous function $ \varphi (x) $ and $ A,B,\varepsilon  > 0 $  one gets 
\begin{flalign}
\nonumber & \int_{ - A}^B {\varphi (x)\left( {\int_0^1 {dt{t^{ix - 1}}{{(1 - t)}^{ - ix - 1}}} } \right)dx}  =  
\\
\nonumber & ~~~~~~~~~~~~~~~~~~~~ = \mathop {\lim }\limits_{\scriptstyle\alpha  \to  + \infty \hfill\atop
\scriptstyle\beta  \to  + \infty \hfill} \int_{ - A}^B {\varphi (x)\left\{ {\frac{1}{{ix}}\left( {\cos (\beta x) - \cos (\alpha x)} \right)} \right\}dx + 2\pi \varphi (0).}
\end{flalign}
The function $ \varphi (x) $ is continuous everywhere except the point $ x = 0, $ and therefore is (absolutely) integrated on any finite interval $ (A,B) $  which does not contain the point $ x = 0, $ For any continuous $ \varphi (x) $ and $ A,B,C > 0, $ such that $ ( - C,C) \subset ( - A,B), $ one has
\begin{eqnarray}
\nonumber \int_{ - A}^B {\varphi (x)\frac{1}{{ix}}\left( {\cos (\beta x) - \cos (\alpha x)} \right)dx}  = \,\,\int_{ - A}^{ - C} { - i} \int_{ - C}^C {\varphi (x)\frac{{\cos (\beta x) - \cos (\alpha x)}}{x}dx}  + \int_C^B .
\end{eqnarray}
Let us consider the second integral in the right hand side separately. Splitting the function $ \varphi (x) $ on the odd and even parts, $ \varphi (x) = {\varphi _{{\rm{even}}}}(x) + {\varphi _{{\rm{odd}}}}(x), $ we obtain
\begin{flalign}
\nonumber & \int_{ - C}^C {\varphi (x)\frac{{\cos (\beta x) - \cos (\alpha x)}}{x}dx}  =
\\
\nonumber & ~~~~~~~~~~ = \,\int_{ - C}^C {{\varphi _{{\rm{even}}}}(x)\frac{{\cos (\beta x) - \cos (\alpha x)}}{x}dx}  + \int_{ - C}^C {{\varphi _{{\rm{odd}}}}(x)\frac{{\cos (\beta x) - \cos (\alpha x)}}{x}dx} . 
\end{flalign}
The first integral in the right hand side is zero, due to the integrand is odd. The second integral can be transformed as:
\begin{eqnarray}
\nonumber ~~~~~ \int_{ - C}^C {{\varphi _{{\rm{odd}}}}(x)\frac{{\cos (\beta x) - \cos (\alpha x)}}{x}dx}  = 2\int_{ - C}^C {\frac{{{\varphi _{{\rm{odd}}}}(x)}}{x}\sin \left( {\frac{{\alpha  - \beta }}{2}x} \right)\sin \left( {\frac{{\alpha  + \beta }}{2}x} \right)dx} .
\end{eqnarray}
The function
\begin{eqnarray}
\nonumber \frac{{{\varphi _{{\rm{odd}}}}(x)}}{x}\sin \left( {\frac{{\alpha  - \beta }}{2}x} \right) \equiv f(x)
\end{eqnarray}
the last expression is continuous and limited in the neighborhood of the point $ x = 0, $ and
\begin{eqnarray}
\nonumber f(0) = {\varphi _{{\rm{odd}}}}(0)\frac{{\alpha  - \beta }}{2} = 0.
\end{eqnarray}
Therefore it is (absolutely) integrated. Hence, according to the Riemann Lemma, one obtains
%\begin{flalign}
%\nonumber & \mathop {\lim }\limits_{\scriptstyle\alpha  \to  + \infty \hfill\atop
%\scriptstyle\beta  \to  + \infty \hfill} \int_{ - C}^C {f(x)\sin \left( {\frac{{\alpha  + \beta }}{2}x} \right)dx}  = 0, \mathop {\lim }\limits_{\scriptstyle\alpha  \to  + \infty \hfill\atop
%\scriptstyle\beta  \to  + \infty \hfill} \int_{ - A}^{ - C} {\varphi (x)\frac{{\cos (\beta x) - \cos (\alpha x)}}{x}dx}  = 0,
%\\
%\nonumber & ~~~~~~~~~~~~~~~~~~~~~~~~ = \mathop {\lim }\limits_{\scriptstyle\alpha  \to  + \infty \hfill\atop
%\scriptstyle\beta  \to  + \infty \hfill} \int_C^B {\varphi (x)\frac{{\cos (\beta x) - \cos (\alpha x)}}{x}dx}  = 0. 
%\end{flalign}
\begin{eqnarray}
\nonumber & \mathop {\lim }\limits_{\scriptstyle\alpha  \to  + \infty \hfill\atop
\scriptstyle\beta  \to  + \infty \hfill} \int_{ - C}^C {f(x)\sin \left( {\frac{{\alpha  + \beta }}{2}x} \right)dx}  = 0, \\ 
\nonumber & \mathop {\lim }\limits_{\scriptstyle\alpha  \to  + \infty \hfill\atop
\scriptstyle\beta  \to  + \infty \hfill} \int_{ - A}^{ - C} {\varphi (x)\dfrac{{\cos (\beta x) - \cos (\alpha x)}}{x}dx}  = 0, \\
\nonumber & \mathop {\lim }\limits_{\scriptstyle\alpha  \to  + \infty \hfill\atop
\scriptstyle\beta  \to  + \infty \hfill} \int_C^B {\varphi (x)\dfrac{{\cos (\beta x) - \cos (\alpha x)}}{x}dx}  = 0. 
\end{eqnarray}

So
\begin{eqnarray}
\nonumber \mathop {\lim }\limits_{\scriptstyle\alpha  \to \infty \hfill\atop
\scriptstyle\beta  \to \infty \hfill} \int_{ - A}^B {\varphi (x)\left\{ {\frac{1}{{ix}}\left( {\cos (\beta x) - \cos (\alpha x)} \right)} \right\}dx = 0.}
\end{eqnarray} \\
{\bf Note 2}
\begin{flalign}
\nonumber & B(\alpha  + ix,\beta  + iy) = \int_0^1 {{t^{\alpha  + ix - 1}}{{(1 - t)}^{\beta  + iy - 1}}dt}  = \\
\nonumber & ~~~~~~~~~~~~~~~~~~ = \int_0^1 {{{\left( {\frac{{1 - t}}{t}} \right)}^{\frac{{\beta  - \alpha  + i(y - x)}}{2}}}{{[t(1 - t)]}^{ - 1 + \frac{{\alpha  + \beta  + i(y + x)}}{2}}}dt} = \\
\nonumber & ~~~~~~~~~~~~~~~~~~~~~~~~~~~~~~ = \int\limits_{ - \infty }^\infty  {\exp \left[ {\xi \frac{{\beta  - \alpha  + i(y - x)}}{2}} \right]{{[2\cosh (\xi /2)]}^{\alpha  + \beta  + i(y + x)}}d\xi ;}
\end{flalign}
\begin{eqnarray}
\nonumber B(\alpha  + ix,\beta  - ix) = \int\limits_{ - \infty }^\infty  {\exp \left[ {\xi \frac{{\beta  - \alpha  - 2ix}}{2}} \right]{{[2\cosh (\xi /2)]}^{\alpha  + \beta }}d\xi ,}
\end{eqnarray} 
\begin{flalign}
\nonumber & B{(ix,iy)_{|y =  - x}} = \mathop {\lim }\limits_{\alpha ,\beta  \to 0 + } B(\alpha  + ix,\beta  - ix) = \\
\nonumber & ~~~~~~~~~~~~~~~~~~~~~~~~~~~~~~~~ = \int\limits_{ - \infty }^\infty  {\exp ( - ix\xi )d\xi  = 2\pi \delta ( - x) = 2\pi \delta (x).}
\end{flalign}

\end{document}